\documentclass{amsart}%
\usepackage{amsfonts}
\usepackage{amsmath}
 \usepackage{amsthm}
\usepackage{hyperref}
\usepackage{amssymb}
\usepackage[lite]{amsrefs}
\usepackage{graphicx}%
\usepackage[margin=1.04in]{geometry}
\newcommand{\Z}{\mathbb{Z}}

\protected\def\verythinspace{%
  \ifmmode
    \mskip0.05\thinmuskip
  \else
    \ifhmode
      \kern0.08334em
    \fi
  \fi
}

\newtheorem{theorem}{Theorem}

\newtheorem*{proposition3}{Proposition 3}

\newtheorem*{theorem4}{Theorem 4}

\newtheorem{example}[theorem]{Example}

\newtheorem{lemma}[theorem]{Lemma}

\newtheorem{proposition}[theorem]{Proposition}
\newtheorem*{remark}{Remark}

\usepackage[foot]{amsaddr}
\begin{document}

\title[On the torsion homology of Oeljeklaus--Toma manifolds]{Studies of the torsion in the homology \\ of Oeljeklaus--Toma manifolds}
\date{\today}
\author{Dung Phuong PHAN$^{1,2}$}
\address{$^1$Laboratoire de math\'ematiques GAATI, Universit\'e de la Polyn\'esie Fran\c{c}aise, BP 6570, 98702 Faa'a, French Polynesia}
\address{$^2$Ho Chi Minh City University of Technology and Education, Vietnam}
\author{Tuan Anh BUI$^{3,4}$}
\address{$^3$Faculty of Mathematics and Computer Science, University of Science, Ho Chi Minh City, Vietnam}
\address{$^4$Vietnam National University, Ho Chi Minh City, Vietnam}
\author{Alexander D. RAHM$^{1,*}$}
\address{$^*$ Alexander.Rahm@upf.pf , ORCID: 0000-0002-5534-2716, https://gaati.org/rahm}

\thanks{This paper is dedicated to Oliver Br\"aunling.}

\subjclass[2010]{
Primary. 53C55: Global differential geometry of Hermitian and Kählerian manifolds.
Secondary. 11F75: Cohomology of arithmetic groups;
11R27: Units and factorization}
\keywords{Oeljeklaus--Toma manifolds; Cohomology of arithmetic groups}

\begin{abstract} This article investigates the torsion homology behaviour in towers of Oeljeklaus--Toma (OT)\ manifolds. This adapts an idea of Silver--Williams from knot theory to OT\ manifolds and extends it to higher degree homology groups.

In the case of surfaces, i.e. Inoue surfaces of type $S^{0}$, the torsion grows exponentially in both $H_{1}$ (as was established by Br\"aunling) and $H_{2}$ (our result) according to a parameter which already plays a role in Inoue's classical paper, and we obtain that the torsion vanishes in all higher degrees.
This motivates our presented machine calculations for OT\ manifolds of one dimension higher.
\end{abstract}

\maketitle

\section{Introduction}

In order to motivate the computations in this paper, let us for a moment look
at an idea from knot theory. Suppose%
\[
K\hookrightarrow S^{3}%
\]
is a tame knot. Then one can put a tubular neighbourhood $T$ around $K$ and
the knot complement is the space%
\[
X_{K}=S^{3}-\left(  \text{interior of }T\right)  \text{.}%
\]
This makes sense not just as a topological space, but indeed gives a compact
connected $3$-manifold. By a standard computation, one always has%
\begin{equation}
H_{1}(X_{K},\mathbb{Z})\simeq\mathbb{Z}\text{,} \label{l1}%
\end{equation}
see for example \cite{GreenbergHarper}*{Remark 18.4}, independently of the choice of
the knot. By the Hurewicz theorem, this means that the quotient of the
fundamental group to its abelianization, i.e. the first homology group of
$X_{K}$, is given by%
\begin{equation}
q:\pi_{1}(X_{K},\ast)\twoheadrightarrow\mathbb{Z}\text{.} \label{l2}%
\end{equation}
This map is almost canonical. Really, in both Equation \ref{l1} and \ref{l2}
one just gets infinite cyclic groups canonically, so for identifying them with
$\mathbb{Z}$ one needs to choose a generator of this infinite cyclic group,
and changing this would change the maps by multiplication with $-1$.

Consider the subgroups%
\[
q^{-1}(n\mathbb{Z})\subseteq\pi_{1}(X_{K},\ast)
\]
for integers $n\geq1$, where $q$ is as in Equation \ref{l2}. These subgroups
are well-defined, since multiplication by $-1$ preserves them. By covering
space theory, each such subgroup corresponds to a finite covering%
\begin{equation}%
\begin{array}
[c]{c}%
X_{K,n}\\
\mid\\
X_{K}%
\end{array}
\label{l4}%
\end{equation}
of degree $n$. Now, one can study the torsion homology of the spaces $X_{K,n}%
$, depending on $n$. Gonz\'{a}les-Acu\~{n}a--Short~\cite{Gonzalez-AcunaShort} and Riley~\cite{Riley} have independently proven in the early 1990s that%
\begin{equation}
\lim_{n\longrightarrow+\infty}\frac{\log\left\vert H_{1}(X_{K,n}%
,\mathbb{Z})_{tor}\right\vert }{n}=\log M(\Delta_{K})\text{,} \label{l3}%
\end{equation}
where $\Delta_{K}$ denotes the Alexander polynomial of the knot, and $M(-)$
the Mahler measure. So, loosely speaking, this shows that the amount of
torsion first homology classes of the spaces $X_{K,n}$ grows in a very
controlled way along $n$. For example, in the case $\left\vert M(\Delta
_{K})\right\vert >1$ one obtains the asymptotic exponential growth%
\[
\left\vert H_{1}(X_{K,n},\mathbb{Z})_{tor}\right\vert \sim M(\Delta_{K}%
)^{n}\qquad\text{as}\qquad n\longrightarrow+\infty\text{.}%
\]
Note that this is just an asymptotic;\ for small $n$ this is a bad heuristic.
The case $\left\vert M(\Delta_{K})\right\vert >1$ is fairly common, so one
sees that the spaces $X_{K,n}$ for large $n$ will have a tremendous amount of
torsion classes in $H_{1}$. More qualitatively, one can also say that once the
torsion homology is not bounded in $n$, it must already grow exponentially,
and if Lehmer's conjecture is true, there is even a lower bound on the minimal
possible exponential growth.

As a brief comment on the literature, we note that many articles instead
discuss the first homology of the branched covering spaces $\widehat{X}_{K,n}$
instead (e.g., \cite{Gonzalez-AcunaShort}), but one just has $H_{1}(X_{K,n}%
,\mathbb{Z})_{tor}\cong H_{1}(\widehat{X}_{K,n},\mathbb{Z})$, so these
considerations are just a different viewpoint, but equivalent, see
\cite{BurdeZieschang}*{Chapter 8}.

In 2002 Daniel Silver and Susan Williams have pointed out in their article~\cite{SilverWilliams} that many of the above considerations do not need to be
restricted to knot complements. Instead, they point out that whenever one has
a connected manifold $X$ with a surjection%
\[
q:\pi_{1}(X,\ast)\twoheadrightarrow\mathbb{Z}%
\]
imitating Equation \ref{l2}, one can run the analogous analysis for any such
$X$. One has the corresponding covering spaces, call them $X_{n}$, and can
study the limit%
\[
\lim_{n\longrightarrow+\infty}\frac{\log\left\vert H_{1}(X_{n},\mathbb{Z}%
)_{tor}\right\vert }{n}\text{.}%
\]
This note is about a special type of complex manifolds where this idea can be implemented:

Concretely, the situation envisioned by Silver--Williams is met for Inoue
surfaces~\cite{Inoue}; and in this note we only look at Inoue surfaces of
so-called type $S^{0}$ (these are the ones discussed in \S 2 loc. cit.). These
are compact non-K\"{a}hler complex surfaces (so instead of dimension 3, we now
move to compact $4$-dimensional real manifolds with rich extra structure.
While the knot complements are frequently hyperbolic, the\ Inoue surfaces
carry at least a locally conformal K\"{a}hler metric. This plays no role here though).

Inoue surfaces $X$ are one of the types appearing in Kodaira's classification of
minimal compact complex surfaces. They occur in the not fully understood Class~VII$_{0}$. They satisfy $H_{1}(X,\mathbb{Q})\simeq\mathbb{Q}$, vaguely
analogous to Equation \ref{l1}, and one can set up a surjection%
\[
q:\pi_{1}(X,\ast)\twoheadrightarrow\mathbb{Z}%
\]
as in Equation \ref{l2}, again canonical up to multiplication with $-1$ in
$\mathbb{Z}$. As pointed out by Br\"aunling and Vuletescu \cite{BraunlingVuletescu}, one obtains%
\begin{equation}
\lim_{n\longrightarrow+\infty}\frac{\log\left\vert H_{1}(X_{n},\mathbb{Z}%
)_{tor}\right\vert }{n}=\log M(f)\text{,} \label{15}%
\end{equation}

where $M(f)$ is the Mahler measure of a polynomial $f$ which in the setting of
Inoue's paper is the minimal polynomial of the matrix he denotes by
\textquotedblleft$M$\textquotedblright\ in \cite{Inoue}*{\S 2}, or, in the
setting of \cite{BraunlingVuletescu} the minimal polynomial of a unit in a certain number field (this translation follows the philosophy of \cite{OeljeklausToma}).

This result is entirely in line with the philosophy laid out by Silver and
Williams in~\cite{SilverWilliams}. This motivates the question to investigate the
remaining torsion homology%
\[
\left\vert H_{r}(X_{n},\mathbb{Z})_{tor}\right\vert
\]
for $r\neq1$.

\begin{theorem}
\label{thm_1}For Inoue surfaces of type $S^{0}$ (resp. Oeljeklaus--Toma
manifolds with $r_{1},r_{2}=1$), the torsion homology growth satisfies%
\[
\log\left\vert H_{r}(X_{n},\mathbb{Z})_{tor}\right\vert \sim\left\{
\begin{array}
[c]{ll}%
n\log M(f) & \text{for }r=1,2\\
0 & \text{for }r=0,3,4
\end{array}
\right.
\]
as $n\longrightarrow+\infty$, where $M(f)$ is the Mahler measure of the
minimal polynomial as described above in the text. This is $>1$, so the
orders of the torsion part of both  $H_{1}$ and $H_{2}$ grow exponentially with $n$.
\end{theorem}

We refer to \S \ref{sect_ProofOfThm1} for the proof. However, Inoue surfaces
of type $S^{0}$ admit a generalization to higher dimensions due to Oeljeklaus
and Toma \cite{OeljeklausToma}. Given a number field $K$ with $r_{1}\geq1$ real
places and $r_{2}=1$ complex places, and any torsion-free finite-index
subgroup%
\[
U\subseteq\mathcal{O}_{K}^{\times,+}%
\]
of the totally positive units (i.e. the units which are positive under any
homomorphism $\mathcal{O}_{K}^{\times}\rightarrow\mathbb{R}$), they attach a
complex manifold%
\[
X(K,U)
\]
which is

\begin{itemize}
\item connected compact of complex dimension $r_{1}+1$,

\item real dimension $2r_{1}+2$,

\item non-K\"{a}hler\ (but carries a locally conformal K\"{a}hler metric),

\item and whose underlying real manifold is a locally symmetric space for a
solvable Lie group.
\end{itemize}

For any such manifold $X(K,U)$ the fundamental group sits in a canonical exact
sequence%
\[
1\longrightarrow {O}_{K}\longrightarrow\pi_{1}(X(K,U),\ast)\overset{q}%
{\longrightarrow}U\longrightarrow1\text{,}%
\]
see \cite{Braunling}*{prop. 6} and also \cite{PartonVuletescu}*{proof of theorem 4.2}. In the special case of $r_{1}=1$, Dirichlet's Unit
Theorem implies that $\mathcal{O}_{K}^{\times,+}\simeq\mathbb{Z}$, canonical
up to multiplication with $-1$, and $X(K,U)$ is merely an Inoue surface. So,
this is just the case as discussed above. However, for $r_{1}=2$, we have%
\[
\mathcal{O}_{K}^{\times,+}\simeq\mathbb{Z}^{2}\text{,}%
\]
canonical up to the action of $\operatorname*{GL}_{2}(\mathbb{Z})$.

A good source of examples stems from fourth roots. Consider the number field%
\[
K=\mathbb{Q}(\sqrt[4]{p})
\]
for a prime $p$. The polynomial $X^{4}-p$ is irreducible by the Eisenstein
criterion, it has two real roots $\pm\sqrt[4]{p}$ and a pair of complex
conjugate roots $\pm i\sqrt[4]{p}$, so that $K$ indeed satisfies $r_{1}=2$ and
$r_{2}=1$. Hence, this is a great source for OT\ threefolds.\\

Let $\sigma_1$ and $\sigma_2$ be the real embeddings and $\sigma_3$ be a complex embedding of $K$ respectively.\\

The ring of integers can be determined by a general result of Gassert
\cite{Gassert}*{Theorem 1.1}. As soon as $p$ satisfies
$ p\neq 1\operatorname{mod}4$
(equivalently, 
$ p^2\neq p\operatorname{mod}4$),
the ring of integers is%
\[
\mathcal{O}_{K}=\mathbb{Z}[\sqrt[4]{p}]
\]
(the theorem is stated in Gassert's paper only as claiming that the ring of
integers is monogenic, but the proof explicitly shows that the index
$[\mathcal{O}_{K}:\mathbb{Z}[\sqrt[4]{p}]]$, using our notation, is one, so
this shows that $\sqrt[4]{p}$ is a concrete generator). The group of units of $O_K$ is
\[
\mathcal{O}_{K}^{\times}\simeq\mathbb{Z}\left\langle u, v%
\right\rangle
\]
for generating units $u, v$. In the example which we study explicitly, namely $p = 2$,
we choose $u=w^2-1$, $v=w+1$, with $w=\sqrt[4]{2}$.\\

Let
\begin{align*}
i = \min\left\{ k \in \mathbb{N}^*: \sigma_1(u^k)>0, \quad  \sigma_2(u^k)>0\right\},\\
j = \min\left\{ k \in \mathbb{N}^*: \sigma_1(v^k)>0, \quad  \sigma_2(v^k)>0\right\},
\end{align*}

then $O_K^{\times,+} =\mathbb{Z} <u^{i},v^{j}>$.\\

For each $m, n \ge 1$, define the monomorphism $\varphi: <u^{im},v^{jn}> \rightarrow Aut(O_K)$, $x \mapsto \varphi(x)$ such that $\varphi(x)(y) = xy$, for all $y \in O_K$. With this monomorphism, the semidirect product $O_K \rtimes_\varphi \mathbb{Z}<u^{im},v^{jn}>$ is defined with the product
\begin{center}
 $(h_1, u^{im_1}v^{jn_1}).(h_2,u^{im_2}v^{jn_2}) = (h_1+u^{im_1}v^{jn_1}h_2,u^{i(m_1+m_2)}v^{j(n_1+n_2)})$,
 \end{center}
for all $h_1, h_2 \in O_K$, $u^{im_1}v^{jn_1}, u^{im_2}v^{jn_2} \in \mathbb{Z}<u^{im},v^{jn}>$.

The action
\begin{align*}
*: (O_K \rtimes_\varphi \mathbb{Z}<u^{im},v^{jn}>) \times (\mathbb{H}^2 \times \mathbb{C}) \rightarrow (\mathbb{H}^2 \times \mathbb{C})\\
((h,k),(z_1,z_2,z_3))\mapsto (\sigma_1(h) + \sigma_1(k)z_1, \sigma_2(h) + \sigma_2(k)z_2,\sigma_3(h) + \sigma_3(k)z_3 )
\end{align*}
is a properly discontinuous action of $O_K \rtimes_\varphi <u^{im},v^{jn}>$ on $\mathbb{H}^2 \times \mathbb{C}$. Then
\begin{align*}
X_{m,n} = X(K, \mathbb{Z}<u^{im},v^{jn}>):=(\mathbb{H}^2 \times \mathbb{C})/(O_K \rtimes_\varphi \mathbb{Z}<u^{im},v^{jn}>)
\end{align*}
is itself an Oeljeklaus--Toma threefold (or ``OT threefold'').\\

In \S \ref{sect_Algorithm}, we describe an algorithm to find the homology of an OT manifold $X_{m,n}$.\\

Next, a detailed algorithm for the case $p=2$ will be described in  \S \ref{sect_Explicitexample}, and a GAP program which implements this algorithm using the HAP package~\cite{Ellis} is distributed online~\cite{Phan}.

\begin{proposition}
Our algorithm produces the following results for the $r$-th degree homology of $X_{m,n}$ when $p=2$:\\
\scriptsize
\begin{tabular}{|l|l|l|l|l|l|l|l|}
\hline
$m$ & $n$ & r=1 & r=2 & r=3 & r=4 & r=5 & r=6 \\
\hline
1 & 1 & $(\Z/2)^2 \oplus \Z^2$  & $ (\Z/2)^4 \oplus (\Z/4)^2 \oplus \Z$ & $ (\Z/2)^4 \oplus (\Z/4)^2$  & $ (\Z/2)^2 \oplus \Z$ & $\Z^2$ & $\Z$ \\
\hline
2 & 1 & $(\Z/2)^2 \oplus \Z^2$  & $(\Z/2)^4 \oplus (\Z/4)^2 \oplus (\Z/3)^2 \oplus \Z$ &  $(\Z/2)^4 \oplus (\Z/4)^2 \oplus (\Z/3)^2 $  & $ (\Z/2)^2 \oplus \Z$  &  $\Z^2$  & $\Z$  \\
\hline
3 &1 &  $(\Z/2)^2 \oplus \Z^2$  & $ (\Z/2)^4 \oplus (\Z/4)^2 \oplus \Z$&  $(\Z/2)^4 \oplus (\Z/4)^2 $ & $ (\Z/2)^2 \oplus \Z$ & $\Z^2$  & $\Z$  \\
\hline
1 & 2 &  $(\Z/2)^2 \oplus \Z^2$ & $ (\Z/2)^4 \oplus (\Z/4)^2 \oplus \Z$& $(\Z/2)^4 \oplus (\Z/4)^2 $ & $ (\Z/2)^2 \oplus \Z$ & $\Z^2$  & $\Z$  \\
\hline
\end{tabular}
\\
\begin{tabular}{|l|l|l|l|}
2 \verythinspace & 2 &  $(\Z/2)^4 \oplus \Z^2$  & $(\Z/2)^6 \oplus (\Z/4)^2 \oplus (\Z/3)^2 \oplus \Z/8 \oplus \Z/16 \oplus \Z$ \\
\hline
3 & 2 &  $(\Z/2)^2 \oplus \Z^2$ & $ (\Z/2)^4 \oplus (\Z/4)^2 \oplus \Z$\\
\hline
1 & 3 &  $(\Z/2)^2 \oplus \Z^2$  & $ (\Z/2)^4 \oplus (\Z/4)^2 \oplus \Z$ \\
\hline
2 & 3 & $(\Z/2)^2 \oplus \Z^2$  & $ (\Z/2)^4 \oplus (\Z/4)^2 \oplus (\Z/9)^2 \oplus (\Z/43)^2 \oplus \Z$\\
\hline
3 & 3 & $(\Z/2)^2 \oplus \Z^2$  & $ (\Z/2)^5 \oplus \Z/4 \oplus \Z/7 \oplus \Z$\\
\hline
\end{tabular}
\normalsize
\end{proposition}

\begin{remark}
 We use Poincar\'e duality and the Universal Coefficient Theorem as in the proof of Theorem 1 to get
 $H^\text{\rm free}_{6-i}(X_{m,n}) \cong H^\text{\rm free}_i(X_{m,n})$ where $i \in \{0,1,2\}$, and
$H^\text{\rm torsion}_{5-i}(X_{m,n}) \cong H^\text{\rm torsion}_i(X_{m,n})$ where $i \in \{0,1,2\}$.
So we can complete the above table in the degrees $r$ in which the authors' machine did run out of memory.
 \end{remark}

Finally, in Section~\ref{sect_Theoreticalresult}, we determine explicitly the first homology of some classes of the manifolds  $X_{m,n}$  in the case $K = \mathbb{Q} (\sqrt[4]{2})$ and obtain the following two propositions.

\begin{proposition}
For Oeljeklaus--Toma manifolds with $r_1=2$, $r_2=1$ when $p=2$, the first homology is of isomorphism type

$H_1(X_{m,1}) \cong \left\{
\begin{matrix}
(\Z/2\Z)^2 \oplus \Z^2 &  \text{for} & 6 \nmid m \\
(\Z/2\Z)^2 \oplus \Z/7\Z \oplus \Z^2 &  \text{for} & 6 \mid m.
\end{matrix}
\right.$
\label{TheProposition}
\end{proposition}

 \begin{theorem} \label{m:odd}
For Oeljeklaus--Toma manifolds with $r_1=2$, $r_2=1$ when $p=2$ and $m$ is odd, the first homology is of isomorphism type $H_1(X_{m,n}) \cong (\Z/2\Z)^2 \oplus (\Z/\beta(m,n)\Z)^2 \oplus \Z^2$ for all $n \in \mathbb{N}^*$, with the integer $\beta(m,n)$ having the following properties.
\begin{itemize}
 \item Let $a_m = \dfrac{1}{2}\left(-1+\sqrt{2}\right)^m + \dfrac{1}{2}\left(-1-\sqrt{2}\right)^m$, which is a negative natural integer. Consider the prime factorization of $-a_m$, which we shall write $-a_m =p_1^{s_1}p_2^{s_2}...p_k^{s_k}$. For each of its primes $p_i$, there is an integer $n(p_i)$ dividing $p_i^2-1$ such that $\beta(m, n)$ is divisible by $p_i$ for $n$ being any multiple of $n(p_i)$. This yields all the factors of $\beta(m,n)$, and hence $\Z/\beta(m, n)\Z = \Z/1\Z = \{0\}$ for $n$ not divisible by any of those periods $n(p_i)$.
\item Let $t_i = \min\left\lbrace \ell \in \mathbb{N}: n \nmid p_i^\ell\cdot n(p_i)\right\rbrace$ and
$v_i = \min\left\lbrace s_i,t_i \right\rbrace$.
\\
Then $\beta(m,n) = p_1^{v_1}p_2^{v_2}...p_k^{v_k}.$
Hence $H_1(X_{m,n})$ is periodic in $n$, with period
\\
$lcm(p_1^{s_1-1}n(p_1),p_2^{s_2-1}n(p_2),...,p_k^{s_k-1}n(p_k))$.
\end{itemize}
\end{theorem}

An explicit description of $\beta(m, n)$ at small values of $m$ is given in Table~\ref{explicit-description} below.
\\
 For the proof of Theorem 4 and Proposition 3, see Section~\ref{sect_Theoreticalresult}.
\\
Let us remark that the periodicity in $n$ of the torsion in the homology described in Theorem~\ref{m:odd} is in analogy to the periodicity in $n$ of the Betti numbers in the framework of Silver and Williams~\cite{Gordon}*{theorem 4.2}.

\begin{proposition}
Our computational results on the isomorphism type of the first homology of some $X_{m,n}$,  in the case $K = \mathbb{Q} (\sqrt[4]{2})$, are the following: \\
\scriptsize
\begin{tabular}{|l|l|l|l|l|l|}
\hline
$m$ & n= 1 & n=2 & n=3 & n=4 & n=5 \\
\hline
1 & $(\Z/2)^2 \oplus \Z^2$ & $(\Z/2)^2 \oplus \Z^2$ &  $(\Z/2)^2 \oplus \Z^2$ &  $(\Z/2)^2 \oplus \Z^2$ &  $(\Z/2)^2 \oplus \Z^2$ \\
\hline
2 &  $(\Z/2)^2 \oplus \Z^2$ & $(\Z/2)^4 \oplus \Z^2$ &  $(\Z/2)^2 \oplus \Z^2$ &  $(\Z/2)^2 \oplus \Z^2$ &  $(\Z/2)^2 \oplus \Z^2$ \\
\hline
3 &  $(\Z/2)^2 \oplus \Z^2$ &  $(\Z/2)^2 \oplus \Z^2$ &  $(\Z/2)^2 \oplus \Z^2$&  $(\Z/2)^2 \oplus \Z^2$ &  $(\Z/2)^2 \oplus \Z^2$\\
\hline
4 &  $(\Z/2)^2 \oplus \Z^2$ & $(\Z/2)^4 \oplus \Z^2$  &  $(\Z/2)^2 \oplus \Z^2$ & $(\Z/4)^2 \oplus (\Z/8)^2 \oplus \Z^2$ & $(\Z/2)^2 \oplus \Z^2$ \\
\hline
5 &  $(\Z/2)^2 \oplus \Z^2$ &  $(\Z/2)^2 \oplus \Z^2$ &  $(\Z/2)^2 \oplus \Z^2$ &  $(\Z/2)^2 \oplus \Z^2$ &  $(\Z/2)^2 \oplus \Z^2$\\
\hline
6 & $(\Z/2)^2 \oplus \Z/7 \oplus \Z^2$ & $(\Z/2)^4 \oplus \Z/7 \oplus \Z^2$  & $(\Z/2)^2 \oplus (\Z/7)^2 \oplus \Z^2$ & $(\Z/2)^4 \oplus \Z/7 \oplus \Z^2$ & $(\Z/2)^2 \oplus \Z/7 \oplus \Z^2$ \\
\hline
7 &  $(\Z/2)^2 \oplus \Z^2$ &  $(\Z/2)^2 \oplus \Z^2$ &  $(\Z/2)^2 \oplus \Z^2$ &  $(\Z/2)^2 \oplus \Z^2$ &  $(\Z/2)^2 \oplus \Z^2$\\
\hline
8 &  $(\Z/2)^2 \oplus \Z^2$ & $(\Z/2)^4 \oplus (\Z/3)^2 \oplus \Z^2$ &  $(\Z/2)^2 \oplus \Z^2$ & $(\Z/3)^4 \oplus (\Z/4)^2 \oplus (\Z/8)^2 \oplus \Z^2$ &  $(\Z/2)^2 \oplus \Z^2$\\
\hline
9 &  $(\Z/2)^2 \oplus \Z^2$ &  $(\Z/2)^2 \oplus \Z^2$ &  $(\Z/2)^2 \oplus \Z^2$ &  $(\Z/2)^2 \oplus \Z^2$ &  $(\Z/2)^2 \oplus \Z^2$\\
\hline
10 &  $(\Z/2)^2 \oplus \Z^2$ & $(\Z/2)^4 \oplus \Z^2$ &  $(\Z/2)^2 \oplus \Z^2$ & $(\Z/2)^4 \oplus \Z^2$ &  $(\Z/2)^2 \oplus \Z^2$\\
\hline
11 &  $(\Z/2)^2 \oplus \Z^2$&  $(\Z/2)^2 \oplus \Z^2$ & $(\Z/2)^2 \oplus \Z^2$ &  $(\Z/2)^2 \oplus \Z^2$ &  $(\Z/2)^2 \oplus \Z^2$\\
\hline
12 & $(\Z/2)^2 \oplus \Z/7 \oplus \Z^2$ & $(\Z/2)^4 \oplus \Z/7 \oplus \Z^2$  & $(\Z/2)^2 \oplus (\Z/7)^2 \oplus \Z^2$ & $ (\Z/4)^2 \oplus (\Z/8)^2 \oplus \Z/7 \oplus \Z^2$ & $(\Z/2)^2 \oplus \Z/7 \oplus \Z^2$ \\
\hline
\end{tabular}
\normalsize
\end{proposition}

We note that recent progress has been made also on other cohomological aspects of OT manifolds, namely on their Dolbeault cohomology~\cites{AgnellaDubickasOtimanStelzig,Kasuya} and on their de~Rham and twisted cohomology~\cite{IstratiOtiman}.

\medskip

\paragraph{\textbf{Acknowledgements.}}
This article is dedicated to Oliver Br\"aunling, for having instigated and prepared the research project documented herein, having established the motivations and described the connection to knot theory, and having taught us the aspects of Oeljeklaus-Toma manifolds which were necessary for our computations. We regret that he persisted not to be named as a co-author, out of modesty. We would like to thank Graham Ellis for support with the HAP package of GAP; Ethan Berkove and Ruben Sanchez-Garcia for their suggestions which led to the elimination of a major flaw in our algorithm, and Nicolina Istrati for a critical lecture of a preliminary manuscript.
We would like to thank the anonymous referee for suggesting some essential improvements of the paper. We acknowledge financial support by the MELODIA project, grant number ANR-20-CE40-0013 of the Agence Nationale de la Recherche.

\begin{table}
\caption{In the setting of Theorem~\ref{m:odd}, we find find the following periods $n(p)$ for the primes $p$ in the factorization of
 $a_m = \dfrac{1}{2}(-1+\sqrt{2})^m + \dfrac{1}{2}(-1-\sqrt{2})^m$,
 all of them dividing $(p-1)\cdot(p+1)$, such that the integer $\beta(m,n)$ of Theorem~\ref{m:odd} is divisible by $p$ for $n$ being any multiple of $n(p)$.}
 \label{explicit-description}
 \begin{tabular}{| p{1cm}| p{5cm} | p{2cm} p{0.5cm} p{4cm}|}
 \hline
 $m$ & $-a_m$ & $\beta(m,n)$ & & \\
 \hline
 $1$ & $1$ & $=1$, & $\forall$ & $n \in \mathbb{N} $.\\
 \hline
 $3$ & $7$ & $ = 7$ & if & $n =24k$; \\
 & & $=1$ & if & $n \neq 24k$. \\
 \hline
 $5$ & $41$ & $=41$ & if & $n = 210k$ \\
 & & $=1$ & if & $n \neq 210k$.\\
 \hline
 $7$ & $239$ & $=239$ & if & $n =119k$ \\
 & & $=1$ & if & $n \neq 119k$.\\
 \hline
 $9$ & $1393= 7\cdot 199$ & $=7$ & if & $n = 24k$ and $n \neq 1800k$ \\
 & & $=1393$ & if & $n = 1800k$ \\
 & & $=1$ & if & $n \neq 24k$.\\
 \hline
 $11$ & $8119= 23\cdot 353$ & $=353$ & if & $ n = 88k$ and $n \neq 264k$ \\
 & & $=8119$ & if & $n = 264k$\\
 & & $= 1$ & if & $n \neq 88k$.\\
 \hline
 $13$ & $47321= 79\cdot 599$ & $=79$ & if & $n = 39k$ and $n \neq 7800k$\\
 & & $=599$ & if & $n = 2600k$ and $n \neq 7800k$\\
 & & $=47321$ & if & $n = 7800k$ \\
 & & $ = 1$ & if & $n \neq 39k$ and $n \neq 2600k$.\\
 & & Note that & &  $47321=79\cdot 599$ and $7800 = lcm(39,2600)$.\\
 \hline
 $15$ & $7\cdot 31^2\cdot 41$ & $=7$ & if & $n = 24k, \ \ \neq 120k$ \\
 & & $=31$ & if & $n =15k, \neq 120k, \neq 465k, \neq 210k$\\
 & & $=7\cdot 31$ & if & $n = 120k, \neq 3720k, \neq 840k$\\
 & & $=31^2$ & if & $n=465k, \neq 3720k, \neq 6510k$ \\
 & & $=31\cdot 41$ & if & $n = 210k, \neq 840k, \neq 6510k$\\
 & & $=7\cdot 31^2$ & if & $n=3720k, \neq 26040k$\\
 & & $=7\cdot 31\cdot 41$ & if & $n = 840k, \neq 26040k$\\
 & & $=31^2\cdot 41$ & if & $n = 6510k, \neq 26040k$\\
 & & $=7\cdot 31^2\cdot 41$ & if & $n =26040k$\\
 & & $=1$ & if & $n \neq 24k, 15k, 210k$.\\
 \hline
 $17$ & $1607521 = 103\cdot 15607$ & \\
 & & $=103$ & if & $n =1768k$\\
 & & $=15607$ & if & $n = 265336k$\\
 & & $=1$ & if & $ n \neq 1768k, \neq 265336k$ \\
 \hline
 $19$ & $9369319$ & $=9369319$ & if & $n=178017080k$ \\
 & & $=1$ & if & $n \neq 178017080k$ \\
 \hline
 \end{tabular}
\end{table}

\section{\label{sect_ProofOfThm1}Proof of Theorem \ref{thm_1}}
Let $X_n=X(K,\mathbb{Z}<v^n>)$. As an OT manifold $X_n$ of dimension 4 is oriented and closed, according to Poincar\'e duality,
\begin{center}
$H^k(X)\cong H_{4-k}(X), 0 \le k \le 4$.
\end{center}
With $k=3$, $H^3(X) \cong H_1(X)$, so the torsion of $H^3(X)$ is isomorphic to the torsion of $H_1(X)$.
On the other hand, by the Universal Coefficient Theorem,
\begin{center}
$H^1(X_n) \cong Torsion(H_0(X_n))\oplus [H_1(X_n)/Torsion(H_1(X_n))]$.
\end{center}
Since $H_0(X_n) \cong \mathbb{Z}$,
\begin{center}
$H_3(X_n) \cong H^1(X_n) \cong [H_1(X_n)/Torsion(H_1(X_n))]$.
\end{center}
This implies
\begin{center}
$Torsion(H_3(X_n)) \cong Torsion(H^1(X_n))=0$.
\end{center}
In the case of $H^2(X_n)$, again with the Universal Coefficient Theorem,
\begin{center}
$H^2(X_n) \cong Torsion(H_1(X_n)) \oplus [H_2(X_n)/Torsion(H_2(X_n))]$.
\end{center}
This implies
\begin{center}
$Torsion(H_2(X_n)) \cong Torsion(H^2(X_n))\cong Torsion(H_1(X_n))$.
\end{center}
Usage of Equation (5), stated above, completes the proof of the theorem.

\newpage

\section{\label{sect_Algorithm}\mbox{ An algorithm computing the homology of $X_{m,n}$ }}
\subsection{A resolution of $\mathbb{Z}$ over $\mathbb{Z}^2$}
Consider the topological space $ \mathbb{R}^2$ and the multiplicative group
$$G = <x,y|\  x \mbox{ and } y \mbox{ commute}>,$$
with the generators acting via
\begin{center}
$(a,b)\mapsto x(a,b):= (a+1,b)$ and $(a,b) \mapsto y(a,b):=(a,b+1)$,
\end{center}
for all $(a,b) \in  \mathbb{R}^2$. Then the group $G$ acts on the space $ \mathbb{R}^2$ by the translation $G \times \mathbb{R}^2: (g,u) \mapsto gu:=g(u)$.\\

With this action, $G$ acts by cellular maps for the following CW-structure on $ \mathbb{R}^2$:
\begin{itemize}
\item 0-cells: all points $ge^0$, with $e^0=(0,0)$, $g \in G$.
\item 1-cells: all intervals 
\begin{align*}
ge^1_1, ge^1_2,
\end{align*}
where $ e^1_1 = \left\{(a,0):0<a<1\right\}$, $ e^1_2 = \left\{(0,b):0<b<1\right\}$, $g \in G$.
\item 2-cells: the squares 
\begin{align*}
ge^2_{12}
\end{align*}
where $e^2_{12} = \left\{(a,b):0<a<1, 0<b<1\right\}$.
\end{itemize}

This CW-structure induces a CW-structure on the quotient space $ \mathbb{R}^2/G$ as following
\begin{itemize}
\item 0-cell: the point $e^0 = (0,0,0,0)$.
\item 1-cells: the intervals
\begin{align*}
e^1_1 = \left\{(a,0):0<a<1\right\}, e^1_2 = \left\{(0,b):0<b<1\right\}
\end{align*}
\item 2-cells: the square
\begin{align*}
e^2_{12} = \left\{ (a,b): 0<a<1, 0<b<1 \right\}
\end{align*}
\end{itemize}
We have $C_*( \mathbb{R}^2) \equiv C_*( \mathbb{R}^2/G)$, where $C_k( \mathbb{R}^2/G)$ is the free abelian group generated by all k-cells in $X/G$.\\

The cellular chain complex
\begin{align*}
0 \rightarrow C_2( \mathbb{R}^2/G)\overset{\partial_2}{\rightarrow}C_1( \mathbb{R}^2/G)\overset{\partial_1}{\rightarrow}C_0( \mathbb{R}^2/G)\overset{\varepsilon}{\rightarrow} \mathbb{Z} \rightarrow 0,
\end{align*}
is a resolution of $\mathbb{Z}$ over $\mathbb{Z}^2$ with  the boundary maps

\begin{align*}
\partial_1(e^1_1) = xe^0 -e^0, \partial_1(e^1_2) = ye^0 -e^0.\\
\partial_2(e^2_{12})= xe^1_2-e^1_2+e^1_1-ye^1_2. 
\end{align*}

and with the contracting homotopies
\begin{itemize}
\item $h_1: C_1 \rightarrow C_2$ with
\begin{align*}
h_1(x^my^ne^1_1) =0; \\
h_1(x^my^ne^1_2)=sign(m)\left(\sum\limits_{i=\overline{m}}^{m-1-\overline{m}} x^i \right)y^ne^2_{12}.
\end{align*}

\item $h_0: C_0 \rightarrow C_1$ with
\begin{align*}
h_0(x^my^ne^0)=sign(m)\left(\sum\limits_{i=\overline{m}}^{m-1-\overline{m}} x^i \right)y^ne^1_1+sign(n)\left(\sum\limits_{j=\overline{n}}^{n-1-\overline{n}} y^j \right)e^1_2.
\end{align*}
\end{itemize}
In the above formulas,
\begin{align*}
\overline{m} = \dfrac{m-|m|}{2}, \quad \overline{n} = \dfrac{n-|n|}{2}.
\end{align*}

Generalizing the above process, we get a resolution of $\mathbb{Z}$ over $\mathbb{Z}^4$.

\subsection{ A resolution of $\mathbb{Z}$ over $\mathbb{Z}^4$}
Consider the topological space $ \mathbb{R}^4$ and the multiplicative group
$$G = <x,\ y,\ z,\ t\ |\ x,\ y,\ z\ \text{and } t\ \text{commute}>,$$
with the generators acting via
\begin{align*}
(a,b,c,d)\mapsto x(a,b,c,d):=(a+1,b,c,d),\\
(a,b,c,d)\mapsto y(a,b,c,d):=(a,b+1,c,d),\\
(a,b,c,d)\mapsto z(a,b,c,d):=(a,b,c+1,d),\\
(a,b,c,d)\mapsto t(a,b,c,d):=(a,b,c,d+1),
\end{align*}
for all $(a,b,c,d) \in \mathbb{R}^4$. Then the group $G$ acts on the space $\mathbb{R}^4$ by the translation $G \times \mathbb{R}^4: (g,u) \mapsto g(u)$.\\
With this action, $G$ acts by cellular maps for the following CW-structure on $\mathbb{R}^4$, and this CW-structure induces a CW-structure on the quotient space $\mathbb{R}^4/G$.
Let $C_k(\mathbb{R}^4/G)$ be the free abelian group generated by all k-cells in $\mathbb{R}^4/G$, then the cellular chain complex
\begin{align*}
0 \rightarrow C_4(\mathbb{R}^4/G)\overset{\partial_4}{\rightarrow}C_3(\mathbb{R}^4/G)\overset{\partial_3}{\rightarrow}C_2(\mathbb{R}^4/G)\overset{\partial_2}{\rightarrow}C_1(\mathbb{R}^4/G)\overset{\partial_1}{\rightarrow}C_0(\mathbb{R}^4/G)\overset{\varepsilon}{\rightarrow} \mathbb{Z} \rightarrow 0,
\end{align*}
is a resolution of $\mathbb{Z}$ over $\mathbb{Z}^4$ with the contracting homotopies
\begin{itemize}
\item $h_3: C_3 \rightarrow C_4$ with
\begin{align*}
h_3(x^my^nz^pt^qe^3_{123}) = h_3(x^my^nz^pt^qe^3_{124})=h_3(x^my^nz^pt^qe^3_{134})=0,\\
h_3(x^my^nz^pt^qe^3_{234})=sign(m)\left(\sum\limits_{i=\overline{m}}^{m-1-\overline{m}} x^i \right)y^nz^pt^qe^4_{1234}.
\end{align*}
\item $h_2: C_2 \rightarrow C_3$, with
\begin{align*}
h_2(x^my^nz^pt^qe^2_{12})=h_2(x^my^nz^pt^qe^2_{13})
=h_2(x^my^nz^pt^qe^2_{14})=0,\\
h_2(x^my^nz^pt^qe^2_{23})=sign(m)\left(\sum\limits_{i=\overline{m}}^{m-1-\overline{m}} x^i \right)y^nz^pt^qe^3_{123},\\
h_2(x^my^nz^pt^qe^2_{24})=sign(m)\left(\sum\limits_{i=\overline{m}}^{m-1-\overline{m}} x^i \right)y^nz^pt^qe^3_{124},\\
h_2(x^my^nz^pt^qe^2_{34})=sign(m)\left(\sum\limits_{i=\overline{m}}^{m-1-\overline{m}} x^i \right)y^nz^pt^qe^3_{134}+sign(n)\left(\sum\limits_{j=\overline{n}}^{n-1-\overline{n}} y^j \right)z^pt^qe^3_{234}.
\end{align*}
\item $h_1: C_1 \rightarrow C_2$, with
\begin{align*}
h_1(x^my^nz^pt^qe^1_1)=0, \\
h_1(x^my^nz^pt^qe^1_2)=sign(m)\left(\sum\limits_{i=\overline{m}}^{m-1-\overline{m}} x^i \right)y^nz^pt^qe^2_{12},\\
h_1(x^my^nz^pt^qe^1_3)=sign(m)\left(\sum\limits_{i=\overline{m}}^{m-1-\overline{m}} x^i \right)y^nz^pt^qe^2_{13}+sign(n)\left(\sum\limits_{j=\overline{n}}^{n-1-\overline{n}} y^j \right)z^pt^qe^2_{23},
\end{align*}
\begin{align*}
h_1(x^my^nz^pt^qe^1_4)=sign(m)\left(\sum\limits_{i=\overline{m}}^{m-1-\overline{m}} x^i \right)y^nz^pt^qe^2_{14}+sign(n)\left(\sum\limits_{j=\overline{n}}^{n-1-\overline{n}} y^j \right)z^pt^qe^2_{24} \\ + sign(p)\left(\sum\limits_{k=\overline{p}}^{p-1-\overline{p}} z^k \right)t^qe^2_{34}.
\end{align*}

\item $h_0: C_0 \rightarrow C_1$ with
\begin{align*}
h_0(x^my^nz^pt^qe^0)=sign(m)\left(\sum\limits_{i=\overline{m}}^{m-1-\overline{m}} x^i \right)y^nz^pt^qe^1_1+sign(n)\left(\sum\limits_{j=\overline{n}}^{n-1-\overline{n}} y^j \right)z^pt^qe^1_2 \\ + sign(p)\left(\sum\limits_{k=\overline{p}}^{p-1-\overline{p}} z^k \right)t^qe^1_3+ sign(q)\left(\sum\limits_{l=\overline{q}}^{q-1-\overline{q}} t^l \right)e^1_4.
\end{align*}
\end{itemize}

In the above formulas,
\begin{align*}
& \overline{m} = \dfrac{m-|m|}{2},\quad \overline{n} = \dfrac{n-|n|}{2}, \quad \overline{p} = \dfrac{p-|p|}{2}, \quad \overline{q} = \dfrac{q-|q|}{2}, \\
& e^1_1=\left\{(a,0,0,0):0<a<1\right\}, \quad e^1_2=\left\{(0,b,0,0):0<b<1\right\}, \\
& e^1_3=\left\{(0,0,c,0):0<c<1\right\},\quad e^1_4=\left\{(0,0,0,d):0<d<1\right\}, \\
& e^2_{ij} = e^1_i \times e^1_j, \quad i \in \left\lbrace 1, 2, 3, 4 \right\rbrace, \quad j \in \left\lbrace i+1,..., 4 \right\rbrace, \\
& e^3_{ijk} = e^1_i \times e^1_j \times e^1_k, \quad i \in \left\lbrace 1, 2, 3, 4 \right\rbrace, \quad j \in \left\lbrace i+1,..., 4 \right\rbrace, \quad k \in \left\lbrace j+1,..., 4 \right\rbrace,\\
& e^4_{1234}=e^1_1 \times e^1_2 \times e^1_3 \times e^1_4.
\end{align*}

\subsection{The isomorphism from $O_K \rtimes_\varphi \mathbb{Z}<u^{im},v^{jn}>$ to $\mathbb{Z}^4 \rtimes_{\overline{\varphi}} \mathbb{Z}^2$}

The multiplicative group $\mathbb{Z}<u^{im},v^{jn}>$ is isomorphic to the additive group~$\mathbb{Z}^2$ by the isomorphism $\alpha: (u^{im})^p(v^{jn})^q \mapsto (p,q)$. The additive group of the ring of integers
\begin{align*}
O_K = \left\{a+bw+cw^2+dw^3 \medspace |\medspace a,b,c,d \in \mathbb{Z}, \ w^4=p\right\}
\end{align*}
of $K$ is isomorphic to $\mathbb{Z}^4$ by the isomorphism $\psi: a+bw+cw^2+dw^3 \mapsto (a,b,c,d)$.\\
Using $\alpha$ and $\psi$, we define $\overline{\varphi}: \mathbb{Z}^2 \rightarrow Aut(\mathbb{Z}^4)$, $(s,t) \mapsto \overline{\varphi}(s,t)$, where
\begin{align*}
\overline{\varphi}(s,t)(y) = \psi(\varphi(\alpha^{-1}(s,t)))\psi^{-1}(y),\quad \text{for all} \quad y \in \mathbb{Z}^4.
\end{align*}
In detail, for $y=(a,b,c,d) \in \mathbb{Z}^4$,
\begin{align*}
\varphi(\alpha^{-1}(s,t))\psi^{-1}(y)
=\left(\begin{matrix}
1 & w & w^2 & w^3
\end{matrix}\right) M^{ims}N^{jnt} y^T
\end{align*}
where $M$ and $N$ are the matrices defined by
\begin{align*}
u \left(\begin{matrix}
1 & w & w^2 & w^3
\end{matrix}\right) = \left(\begin{matrix}
1 & w & w^2 & w^3
\end{matrix}\right)M,\\
v \left(\begin{matrix}
1 & w & w^2 & w^3
\end{matrix}\right) = \left(\begin{matrix}
1 & w & w^2 & w^3
\end{matrix}\right)N
\end{align*}
respectively.
Finally, $\overline{\varphi}(s,t)(y) = \psi(\varphi(\alpha^{-1}(s,t)))\psi^{-1}(y)=(A_{s,t}, B_{s,t}, C_{s,t}, D_{s,t})$ with
\begin{align*}
\left(\begin{matrix}
A_{s,t}\\
B_{s,t}\\
C_{s,t}\\
D_{s,t}
\end{matrix} \right)
=
M^{ims}
N^{jnt}
y^T.
\end{align*}

This homomorphism helps us to form the semidirect product $\mathbb{Z}^4 \rtimes_{\overline{\varphi}} \mathbb{Z}^2$ which is isomorphic to
\\
$O_K \rtimes_\varphi \mathbb{Z}<u^{im},v^{jn}>$.
\\
Since two groups which are isomorphic have the same homology groups, to find the homology groups of $O_K \rtimes_\varphi \mathbb{Z}<u^{im},v^{jn}>$, we just need to find the homology groups of $\mathbb{Z}^4 \rtimes_{\overline{\varphi}} \mathbb{Z}^2$. The following is the algorithm to find the latter ones.

\subsection{An algorithm for computing the homology of $\mathbb{Z}^4 \rtimes_{\overline{\varphi}} \mathbb{Z}^2$}
\begin{itemize}
\item[\textbf{Step 1.}] Construct $\overline{\varphi}: \mathbb{Z}^2 \rightarrow Aut(\mathbb{Z}^4)$, $(s,t) \mapsto \overline{\varphi}(s,t)$, with $\overline{\varphi}(s,t):(a,b,c,d) \mapsto (A_{s,t},B_{s,t},C_{s,t}, D_{s,t})$, where $A_{s,t}, B_{s,t}, C_{s,t}, D_{s,t}$ are defined by the following formula:
\begin{align*}
\left(\begin{matrix}
A_{s,t}\\
B_{s,t}\\
C_{s,t}\\
D_{s,t}
\end{matrix} \right)
=
M^{ims}
N^{jnt}
\left(\begin{matrix}
a\\
b\\
c\\
d
\end{matrix}\right).
\end{align*}
Then construct the semidirect product $E = \mathbb{Z}^4 \rtimes_{\overline{\varphi}} \mathbb{Z}^2$ with this $\overline{\varphi}$.
\item[\textbf{Step 2.}] The projection $p: E \rightarrow \mathbb{Z}^2$, $(a,b,c,d,e,f) \mapsto (e,f)$, \mbox{for all $(a,b,c,d,e,f) \in E$} has $\text{ker}(p) \cong \mathbb{Z}^4$. So $E$ is a group extension of its normal subgroup $\mathbb{Z}^4$ and the quotient group $\mathbb{Z}^2$. Let $R$ and $S$ be the free resolution of $\mathbb{Z}$ over $\mathbb{Z}^4$ and $\mathbb{Z}^2$ respectively. The twisted tensor product of these resolutions, which we obtain using C.T.C Wall's method~\cite{Wall}, gives us a free resolution for the semidirect product $E$.
\item[\textbf{Step 3.}] Compute the desired homology groups using the above resolution.
\end{itemize}

\section{\label{sect_Explicitexample} An explicit example when $p=2$}
In this case, the field $K$ is $\mathbb{Q}(\sqrt[4]{2})$. Two generators of the free part of $O_K^{\times}$ are $u=w^2-1$, $v=w+1$, with $w=\sqrt[4]{2}$, and $O_K^{\times}=\mathbb{Z}<u,v>$.\\
Since $\sigma_1(u)>0$, $\sigma_2(u)>0$, $\sigma_1(v)>0$, $\sigma_2(v)<0$, $O_K^{\times,+} = \mathbb{Z}<u,v^2>$.\\

We observe that \begin{align*}
u \left(\begin{matrix}
1 & w & w^2 & w^3
\end{matrix}\right) = \left(\begin{matrix}
1 & w & w^2 & w^3
\end{matrix}\right) 
M,
\end{align*}
\begin{align*}
v \left(\begin{matrix}
1 & w & w^2 & w^3
\end{matrix}\right) =  \left(\begin{matrix}
1 & w & w^2 & w^3
\end{matrix}\right) 
N,
\end{align*}
with 
\begin{align*}
M = \left(\begin{matrix}
-1 & 0 & 2 & 0\\
0 & -1 & 0 & 2\\
1 & 0 & -1 & 0\\
0 & 1 & 0 & -1
\end{matrix}\right), \quad N = \left(\begin{matrix}
1 & 0 & 0 & 2\\
1 & 1 & 0 & 0 \\
0 & 1 & 1 & 0\\
0 & 0 & 1 & 1
\end{matrix}\right).
\end{align*}
So, $O_K \rtimes_\varphi \mathbb{Z}<u^m,v^{2n}>$ is isomorphic to $\mathbb{Z}^4 \rtimes_{\overline{\varphi}} \mathbb{Z}^2$ by the isomorphism 
\\
$\overline{\varphi}(s,t)(y) = \psi(\varphi(\alpha^{-1}(s,t)))\psi^{-1}(y)=(A_{s,t}, B_{s,t}, C_{s,t}, D_{s,t})$ with
\begin{equation}
\left(\begin{matrix}
A_{s,t}\\
B_{s,t}\\
C_{s,t}\\
D_{s,t}
\end{matrix} \right)
=
M^{ms}
N^{2nt}
\left(\begin{matrix}
a\\
b\\
c\\
d
\end{matrix}\right).
\label{productInTermsOfMandN}
\end{equation}

\section{\label{sect_Theoreticalresult}Proof of Theorem 4 and Proposition 3}
Let us now provide the proofs of Proposition 3 and Theorem 4 stated in the Introduction. We keep the notations of the previous sections.
\begin{proposition3}
For Oeljeklaus--Toma manifolds with $r_1=2$, $r_2=1$ when $p=2$, the first homology is of isomorphism type

$H_1(X_{m,1}) \cong \left\{
\begin{matrix}
\Z/2\Z \oplus \Z/2\Z \oplus (\Z)^2 &  \text{for} & 6 \nmid m \\
\Z/2\Z \oplus \Z/14\Z \oplus (\Z)^2 &  \text{for} & 6 \mid m.
\end{matrix}
\right.$
\end{proposition3}

\begin{proof}
Consider the following presentations of groups by generators and relations:
\begin{align*}
& < X \medspace | \medspace R > := \ < f_1, f_2,f_3,f_4 \medspace | \medspace f_if_j = f_jf_i,\ i, j \in \left\lbrace 1, 2, 3, 4 \right\rbrace> \ \cong \mathbb{Z}^4, \\
& <Y \medspace | \medspace S> := \ <u_1,u_2 \medspace | \medspace u_1u_2=u_2u_1> \ \cong \mathbb{Z}^2.
\end{align*}
Then we obtain the group
\begin{align*}
E = \mathbb{Z}^4 \rtimes_{\overline{\varphi}} \mathbb{Z}^2= <X \cup Y \medspace | \medspace R \cup S \cup \left\lbrace yxy^{-1} =\overline{\varphi}(y)(x), \forall x \in X, \forall y \in Y \right\rbrace>.
\end{align*}
where with $y = (s,t)$, we apply Equation~\ref{productInTermsOfMandN}:
\begin{align*}
\overline{\varphi}(y)(x) = \left(
  \begin{matrix}
  -1 & 0 & 2 & 0 \\
  0 & -1 & 0 & 2 \\
  1 & 0 & -1 & 0 \\
  0 & 1 & 0 & -1
  \end{matrix} \right)^{ms}
  \left(\begin{matrix}
  1 & 0 & 0 & 2 \\
  1 & 1 & 0 & 0\\
  0 & 1 &1 & 0\\
  0 & 0 & 1 & 1
  \end{matrix}\right)^{2t} x.
\end{align*}
We can compute $H_1(E) \cong E/[E,E]$ as the commutator factor group. To find the commutator $[E,E]$, we consider the system of equations 
  \begin{center}
  $u_if_ju_i^{-1} = \varphi (u_i)(f_j) = f_1^{a_{ij}}f_2^{b_{ij}}f_3^{c_{ij}}f_4^{d_{ij}}, \quad i \in \left\lbrace 1, 2 \right\rbrace, \quad j \in \left\lbrace 1, 2, 3, 4 \right\rbrace$,
  \end{center}
  from which we can imply a system of equations of equivalence classes
  \begin{center}
  $ \overline{f_1}^{a_{ij}-e_{1j}} \overline{f_2}^{b_{ij}-e_{2j}} \overline{f_3}^{c_{ij}-e_{3j}} \overline{f_4}^{d_{ij}-e_{4j}} = \overline{1}, \quad i \in \left\lbrace 1, 2 \right\rbrace, \quad j \in \left\lbrace 1, 2, 3, 4 \right\rbrace$
  \end{center}
  where $e_{kj} = \left\{
  \begin{matrix}
  1 & \text{if} & k = j \\
  0 & \text{if} & k \neq j.
  \end{matrix}\right.$ \\

  \bigskip

  Picking up the power of $\left\{\overline{f_i} \right\}$ in each equation, we obtain a matrix of size $8 \times 4$. In the case $m$ is arbitrary and $n =1$, the matrix has the form
  \begin{align*}
  \left(
  \begin{matrix}
  a_m -1 & 0 & b_m & 0 \\
  0 & a_m -1 & 0 & b_m \\
  2b_m & 0 & a_m-1 & 0 \\
  0 & 2b_m & 0 & a_m -1 \\
  0 & 2 & 1 & 0 \\
  0 & 0 & 2 & 1\\
  2 & 0 & 0 & 2\\
  4 & 2 & 0 & 0
  \end{matrix}
  \right)
  \end{align*}
  where
\\${}\hfill
   a_m = \dfrac{1}{2}(-1+\sqrt{2})^m + \dfrac{1}{2}(-1-\sqrt{2})^m,\hfill
   b_m = \dfrac{a_{m+1} +a_m}{2}.
 \hfill$\\

 \bigskip

  Applying row transformations, we change this matrix into the following matrix (*):
  \begin{align*}
  \left(
  \begin{matrix}
  2 & 0 & 0 & 2 \\
  0 & 2 & 1 & 0 \\
  0 & 0 & 1 & 4\\
  0 & 0 & 0 & 7\\
  0 & 0 & 0 & a_m+4b_m-1 \\
  0 & 0 & 0 & 2(a_m-1)+b_m
  \end{matrix}
  \right) \quad (*)
  \end{align*}
  Now, when $m = 6k, k \in \mathbb{N}^*$, then $gcd(7,a_m+4b_m-1,2(a_m-1)+b_m)=7$, and the matrix (*) is row equivalent to the matrix
  \begin{align*}
   \left(
  \begin{matrix}
  2 & 0 & 0 & 2 \\
  0 & 2 & 1 & 0 \\
  0 & 0 & 1 & 4\\
  0 & 0 & 0 & 7\\
  0 & 0 & 0 & 0\\
  0 & 0 & 0 & 0\\
  \end{matrix}
  \right).
  \end{align*}
 In this case, $H_1(X_{m,1}) \cong (\Z/2\Z)^2 \oplus \Z/7\Z \oplus (\Z)^2$.\\
 On the other hand, when $m \neq 6k, k \in \mathbb{N}^*$, then
 \begin{align*}
 gcd(7,a_m+4bm-1,2(a_m-1)+b_m)=1,
 \end{align*}
 and by applying some row and column transformations, the matrix (*) turns into
 
  \begin{align*}
   \left(
  \begin{matrix}
  2 & 0 & 0 & 2 \\
  0 & 2 & 1 & 0 \\
  0 & 0 & 1 & 4\\
  0 & 0 & 0 & 1\\
  0 & 0 & 0 & 0\\
  0 & 0 & 0 & 0\\
  \end{matrix}
  \right).
  \end{align*}
  In this case, $H_1(X_{m,1})\cong  (\Z/2\Z)^2 \oplus (\Z)^2$.
\end{proof}

\newpage
For the proof of Theorem~\ref{m:odd}, we need the following two lemmata.

\begin{lemma} \label{Lemma 1} Suppose that $\gamma$ is a fourth root of $2$ and $p$ is a prime. Then $(1+\gamma)^{p^2-1} \equiv 1$ mod $p$.
\end{lemma}

\begin{proof}
Let $\mathbb{F}_p$ be the prime field of characteristic $p$ and $\overline{\mathbb{F}_p}$ be the algebraic closure of $\mathbb{F}_p$. Since $\gamma$ is a root of the equation $x^4-2=0$ in $\overline{\mathbb{F}_p}[x]$, then $\gamma \in \overline{\mathbb{F}_p}$.\\
By the Frobenius homomorphism in the field $\overline{\mathbb{F}_p}$,
\begin{align}\label{lema.1}
(x+y)^p = x^p + y^p
\end{align}
for all $x, y \in \overline{\mathbb{F}_p}$.\\
Apply \ref{lema.1} with $x=1$ and $y=\gamma$,
\begin{align}
(1+\gamma)^p = 1 + \gamma^p \\
\Rightarrow (1+\gamma)^{p^2} = (1+\gamma^p)^p = 1+ \gamma^{p^2}.
\end{align}
Let $A = \left\lbrace 1 + \gamma \ |\ \gamma \ \text{ is a root of equation} \ x^4-2=0 \right\rbrace $. Using Fermat's little theorem ($a^{p-1} \equiv 1 \mod p$ if $a$ is not divisible by $p$), from $\gamma$ being a root of $x^4-2=0$, we can deduce that $\gamma^{p^2}$ is also a root of this equation, as follows: $\gamma^{p^2-1} = (\gamma^{p+1})^{p-1} \equiv 1 \mod p$, therefore $0 \equiv (\gamma^{p^2} - \gamma)^4 \equiv (\gamma^{p^2})^4 - \gamma^4  \mod p$.
 This means that $A$ is invariant under the power $2$ of the Frobenius homomorphism. And we can conclude that $A$ belongs to a finite field of degree dividing 2 over $\mathbb{F}_p$. So,
\begin{align}
(1+\gamma)^{p^2-1} =1.
\end{align}
Finally, in the ring of integers,
\begin{align}
(1+\gamma)^{p^2-1} \equiv 1 \ \text{mod} \ p.
\end{align}
\end{proof}

\begin{lemma} \label{Lemma 2} Let $p$ be a prime, then with $n=p^2-1$, and $C_{n}^{i}=\dfrac{n!}{i!(n-i)!}$ the number of $i$-combinations in the set of $n$ elements,
\begin{align*}
-1+ e_n = -1+ \sum\limits_{i=0}^{(p^2-1)/2} 2^i C_{2n}^{4i}, \\
f_n = \sum\limits_{i=0}^{(p^2-1)/2} 2^i C_{2n}^{4i+1}, \\
g_n = \sum\limits_{i=0}^{(p^2-1)/2} 2^i C_{2n}^{4i+2}, \\
h_n = \sum\limits_{i=0}^{(p^2-1)/2} 2^i C_{2n}^{4i+3}
\end{align*}
are all divisible by $p$.
\end{lemma}

\begin{proof}
First, consider the fourth root $i$ of $1$, we have
\begin{align*}
\dfrac{1}{4}\sum\limits_{j=0}^3  i^{kj} =
\left\{
\begin{matrix}
1 & \text{if} & k \equiv 0 \ \text{mod} \ 4 \\
0 & \text{if} & \text{otherwise}.
\end{matrix} \right.
\end{align*}
Now,
\begin{align*}
& e_n = \sum\limits_{i=0}^{(p^2-1)/2} 2^i C_{2n}^{4i} = \sum\limits_{i=0}^{(p^2-1)/2} (2^{1/4})^{4i} C_{2n}^{4i} = \sum\limits_{k=0,\ k\equiv 0 (4)}^{2(p^2-1)} (2^{1/4})^k C_{2n}^k \\
& e_n = \sum\limits_{k=0}^{2(p^2-1)} \dfrac{1}{4}\sum\limits_{j=0}^3  i^{kj} (2^{1/4})^k C_{2n}^k  =\dfrac{1}{4}\sum\limits_{j=0}^3  \sum\limits_{k=0}^{2(p^2-1)} (i^j2^{1/4})^k C_{2n}^k \\
& = \dfrac{1}{4}\sum\limits_{j=0}^3 (1+i^j2^{1/4})^{2(p^2-1)}.
\end{align*}
As $i^j2^{1/4}$ is a fourth root of $2$, by Lemma~\ref{Lemma 1},
\begin{align*}
(1+ i^j2^{1/4})^{p^2-1} \equiv 1 \ \text{mod} \ p.
\end{align*}
This yields that
\begin{align*}
e_n= \dfrac{1}{4}\sum\limits_{j=0}^3 (1+i^j2^{1/4})^{2(p^2-1)} \equiv \dfrac{1}{4} \sum\limits_{j=0}^3 1^2 \\
\equiv 1 \ \text{mod} \ p.
\end{align*}
or $e_n -1 \equiv 0 \ \text{mod} \ p$.\\
Next, we can easily see that
\begin{align*}
& f_n = \sum\limits_{i=0}^{(p^2-1)/2} 2^i C_{2n}^{4i+1} = \sum\limits_{k=0,\ k\equiv 1 (4)}^{2(p^2-1)} (2^{1/4})^{k-1} C_{2n}^k \\
& f_n = \sum\limits_{k=0}^{2(p^2-1)} \dfrac{1}{4}\sum\limits_{j=0}^3  i^{j(k+3)} (2^{1/4})^{k-1} C_{2n}^k  =\dfrac{1}{4\cdot 2^{1/4}}\sum\limits_{j=0}^3 i^{3j} \sum\limits_{k=0}^{2(p^2-1)} (i^j2^{1/4})^k C_{2n}^k \\
& = \dfrac{1}{4\cdot 2^{1/4}}\sum\limits_{j=0}^3 i^{3j}(1+i^j2^{1/4})^{2(p^2-1)} \equiv \dfrac{1}{4\cdot 2^{1/4}}\sum\limits_{j=0}^3 i^{3j} \equiv 0 \ \text{mod} \ p.
\end{align*}
and
\begin{align*}
& g_n = \sum\limits_{i=0}^{(p^2-1)/2} 2^i C_{2n}^{4i+2} = \sum\limits_{k=0,\ k\equiv 2 (4)}^{2(p^2-1)} (2^{1/4})^{k-2} C_{2n}^k \\
& g_n = \sum\limits_{k=0}^{2(p^2-1)} \dfrac{1}{4}\sum\limits_{j=0}^3  i^{j(k+2)} (2^{1/4})^{k-2} C_{2n}^k  =\dfrac{1}{4\cdot 2^{1/2}}\sum\limits_{j=0}^3 i^{2j} \sum\limits_{k=0}^{2(p^2-1)} (i^j2^{1/4})^k C_{2n}^k \\
& = \dfrac{1}{4\cdot 2^{1/2}}\sum\limits_{j=0}^3 i^{2j}(1+i^j2^{1/4})^{2(p^2-1)} \equiv \dfrac{1}{4\cdot 2^{1/2}}\sum\limits_{j=0}^3 i^{2j} \equiv 0 \ \text{mod} \ p.
\end{align*}
and finally
\begin{align*}
& h_n = \sum\limits_{i=0}^{(p^2-1)/2} 2^i C_{2n}^{4i+3} = \sum\limits_{k=0,\ k\equiv 3 (4)}^{2(p^2-1)} (2^{1/4})^{k-3} C_{2n}^k \\
& h_n = \sum\limits_{k=0}^{2(p^2-1)} \dfrac{1}{4}\sum\limits_{j=0}^3  i^{j(k+1)} (2^{1/4})^{k-3} C_{2n}^k  =\dfrac{1}{4\cdot 2^{3/4}}\sum\limits_{j=0}^3 i^{j} \sum\limits_{k=0}^{2(p^2-1)} (i^j2^{1/4})^k C_{2n}^k \\
& = \dfrac{1}{4\cdot 2^{3/4}}\sum\limits_{j=0}^3 i^{j}(1+i^j2^{1/4})^{2(p^2-1)} \equiv \dfrac{1}{4\cdot 2^{3/4}}\sum\limits_{j=0}^3 i^{j} \equiv 0 \ \text{mod} \ p.
\end{align*}
\end{proof}

 \begin{theorem4}
For Oeljeklaus--Toma manifolds with $r_1=2$, $r_2=1$ when $p=2$ and $m$ is odd, the first homology is of isomorphism type $H_1(X_{m,n}) \cong (\Z/2\Z)^2 \oplus (\Z/\beta(m,n)\Z)^2 \oplus \Z^2$ for all $n \in \mathbb{N}^*$, with the integer $\beta(m,n)$ having the following properties.
\begin{itemize}
 \item Let $a_m = \dfrac{1}{2}\left(-1+\sqrt{2}\right)^m + \dfrac{1}{2}\left(-1-\sqrt{2}\right)^m$, which is a negative natural integer. Consider the prime factorization of $-a_m$, which we shall write $-a_m =p_1^{s_1}p_2^{s_2}...p_k^{s_k}$. For each of its primes $p_i$, there is an integer $n(p_i)$ dividing $p_i^2-1$ such that $\beta(m, n)$ is divisible by $p_i$ for $n$ being any multiple of $n(p_i)$. This yields all the factors of $\beta(m,n)$, and hence $\Z/\beta(m, n)\Z = \Z/1\Z = \{0\}$ for $n$ not divisible by any of those periods $n(p_i)$.
\item Let $t_i = \min\left\lbrace \ell \in \mathbb{N}: n \nmid p_i^\ell\cdot n(p_i)\right\rbrace$ and
$v_i = \min\left\lbrace s_i,t_i \right\rbrace$.
\\
Then $\beta(m,n) = p_1^{v_1}p_2^{v_2}...p_k^{v_k}.$
Hence $H_1(X_{m,n})$ is periodic in $n$, with period
\\
$lcm(p_1^{s_1-1}n(p_1),p_2^{s_2-1}n(p_2),...,p_k^{s_k-1}n(p_k))$.
\end{itemize}
\end{theorem4}

\begin{proof}
Applying an analogous procedure as in the proof of Proposition \ref{TheProposition} for the case $m$ is odd and $n$ is arbitrary, we obtain the matrix of powers of $\left\{\overline{f_i}\right\}$ of the form
 \begin{align*}
  \left(
  \begin{matrix}
  a_m -1 & 0 & b_m & 0 \\
  0 & a_m -1 & 0 & b_m \\
  2b_m & 0 & a_m-1 & 0 \\
  0 & 2b_m & 0 & a_m -1 \\
 e_n -1 & f_n & g_n & h_n \\
 2h_n & e_n -1 & f_n & g_n \\
 2g_n & 2h_n & e_n -1 & f_n \\
 2f_n & 2g_n & 2h_n & e_n -1
  \end{matrix}
  \right)
\end{align*}
 where with $k = [n/2]$ the integer part of $n/2$, and $C_{n}^{i}=\dfrac{n!}{i!(n-i)!}$ the number of $i$-combinations in the set of $n$ elements,
  \begin{align*}
  e_n = \sum\limits_{i =0}^k 2^i C_{2n}^{4i}, \ f_n = \sum\limits_{i =0}^k 2^i C_{2n}^{4i+1}, \ g_n = \sum\limits_{i =0}^k 2^i C_{2n}^{4i+2}, \ h_n = \sum\limits_{i =0}^k 2^i C_{2n}^{4i+3}.
  \end{align*}
  Here, we make the convention that if $i > n$, then $C_{n}^{i} = 0$.
 By applying some row operations $ \alpha r_i + r_j \rightarrow r_j$ and two appropriate column operations, this matrix turns into
 \begin{align*}
  \left(
  \begin{matrix}
  2 & 0 & 0 & 0 \\
 0 & 2 & 0 & 0 \\
 0 & 0 & -a_m & 0 \\
 0 & 0 & 0 & -a_m \\ 
0 & 0 & (1-e_n)(a_m+b_m)-g_n & -f_n(a_m+b_m)-h_n \\
 0 & 0 & -2h_n(a_m+b_m)-f_n & (1-e_n)(a_m+b_m)-g_n \\
 0 & 0 & -2g_n(a_m+b_m)+1-e_n & -2h_n(a_m+b_m)-f_n \\
 0 & 0 & -2f_n(a_m+b_m)-2h_n & -2g_n(a_m+b_m)+1-e_n
  \end{matrix}
  \right):=T^{m,n}.
\end{align*}
 Let 
$$ \alpha(m,n) =\ gcd(T^{m,n}_{33}, T^{m,n}_{53},T^{m,n}_{63}, T^{m,n}_{73},T^{m,n}_{83}),$$ \begin{equation} \label{beta:def}
 \beta(m,n) =\ gcd(T^{m,n}_{44}, T^{m,n}_{54},T^{m,n}_{64}, T^{m,n}_{74},T^{m,n}_{84}).
 \end{equation}
 Then the homology is of isomorphism type
 $$H_1(X_{m,n}) \cong (\Z/2\Z)^2 \oplus \Z/(\alpha(m,n) \Z) \oplus \Z/(\beta(m,n) \Z) \oplus \Z^2.$$
 As $m$ is odd, $a_m$ is odd, so comparing the entries of the matrix $T^{m,n}$, we see that $\beta(m,n) = \alpha(m,n)$.

Consider the prime factorization of $-a_m$, which we shall write $-a_m =p_1^{s_1}p_2^{s_2}...p_k^{s_k}$.
For each prime $p_i$, by Lemma~\ref{Lemma 2},
\begin{align*}
p_i \mid e_n-1, f_n, g_n, h_n,
\end{align*}
with $n=p_i^2-1$. This implies that $p_i \mid \beta(m,p_i^2-1)$. Now we can see that there exists at least one number $n$ so that $\beta(m,n)$ is divisible by the prime $p_i$. Choose $n(p_i)$ to be the smallest such $n$.

\begin{example}
 In the case $m =3$ of Theorem~\ref{m:odd}, we have $\beta(3,n)=\begin{cases}
                                                  1, & \text{for } 24 \nmid n, \\
7, & \text{for }  24 \mid n.
                                                 \end{cases}
$
\end{example}
We have $a_m = -7$, $b_m = 5$, $a_m + b_m =-2$ and $\alpha(3,n) \in \left\lbrace 1, 7 \right\rbrace$.\\
 By case-by-case calculations, we can see that $\alpha(3,n) = 1$ for $n$ in $\left\lbrace 1,2,...,23 \right\rbrace$ and $\alpha(3,24) = 7$.\\
 With a fixed value $n$, for $i \in \mathbb{N}^*$,
 \begin{align*}
 \left( \begin{matrix}
 e_{n+i} \\ f_{n+i} \\ g_{n+i} \\ h_{n+i}
\end{matrix} \right) = (N^2)^i
\left(
\begin{matrix}
e_n \\ f_n \\ g_n \\ h_n
\end{matrix}  \right) =
\left( \begin{matrix}
e_i & 2h_i & 2g_i & 2f_i \\
f_i & e_i & 2h_i & 2g_i \\
g_i & f_i & e_i & 2h_i \\
h_i & g_i & f_i & e_i
\end{matrix} \right) \left(
\begin{matrix}
e_n \\ f_n \\ g_n \\ h_n
\end{matrix}  \right).
 \end{align*}
 Applying this expression helps us to obtain
 \begin{equation}
\left( \begin{matrix}
 T^{3,n+i}_{53} \\ T^{3,n+i}_{63} \\ T^{3,n+i}_{73} \\ T^{3,n+i}_{83}
 \end{matrix} \right) = \left( \begin{matrix}
e_i & 2h_i & 2g_i & 2f_i \\
f_i & e_i & 2h_i & 2g_i \\
g_i & f_i & e_i & 2h_i \\
h_i & g_i & f_i & e_i
\end{matrix} \right)^T \left( \begin{matrix}
 T^{3,n}_{53} \\ T^{3,n}_{63} \\ T^{3,n}_{73} \\ T^{3,n}_{83}
 \end{matrix} \right) + \left( \begin{matrix}
 T^{3,i}_{53} \\ T^{3,i}_{63} \\ T^{3,i}_{73} \\ T^{3,i}_{83}
 \end{matrix} \right).
 \label{matrixequality}
 \end{equation}
 This matrix equality yields that if $\alpha(3,n) = 7$, then
 \begin{equation}
 \left( \begin{matrix}
 T^{3,n+i}_{53} \\ T^{3,n+i}_{63} \\ T^{3,n+i}_{73} \\ T^{3,n+i}_{83}
 \end{matrix} \right) \equiv \left( \begin{matrix}
 T^{3,i}_{53} \\ T^{3,i}_{63} \\ T^{3,i}_{73} \\ T^{3,i}_{83}
 \end{matrix} \right) \text{mod} \ 7,\ \text{for} \ i \in \mathbb{N}^*.
 \label{samemod}
 \end{equation}
 Choose $i = 24$,
 \begin{align*}
 \left( \begin{matrix}
 T^{3,n+24}_{53} \\ T^{3,n+24}_{63} \\ T^{3,n+24}_{73} \\ T^{3,n+24}_{83}
 \end{matrix} \right) \equiv \left( \begin{matrix}
 T^{3,24}_{53} \\ T^{3,24}_{63} \\ T^{3,24}_{73} \\ T^{3,24}_{83}
 \end{matrix} \right) \equiv
 \left( \begin{matrix}
 0 \\ 0 \\ 0 \\ 0 \\
\end{matrix}  \right) \text{mod} \ 7.
 \end{align*}
 This means that $\alpha(3,n+24) = \alpha(3,24) = 7$.\\
 Next, if $\alpha(3,n) = 1$, then $n = 24l + i$, with $l \in \mathbb{N}$ and $ i \in \left\lbrace 1,2,...,23 \right\rbrace$. As $\alpha(3,24l+24) = 7$, Equation~\ref{samemod} shows that
 \begin{equation}
 \left( \begin{matrix}
 T^{3,24l+i+24}_{53} \\ T^{3,24l+i+24}_{63} \\ T^{3,24l+i+24}_{73} \\ T^{3,24l+i+24}_{83}
 \end{matrix} \right) \equiv \left( \begin{matrix}
 T^{3,i}_{53} \\ T^{3,i}_{63} \\ T^{3,i}_{73} \\ T^{3,i}_{83}
 \end{matrix} \right) \text{mod} \ 7,\ \text{for} \ i \in \mathbb{N}^*.
 \label{samemodn24l}
 \end{equation}
 This means that $\alpha(3,n+24) = \alpha(3,24l+i+24)=\alpha(3,i)$. As $\alpha(3,i) = 1$ for $i \in \left\lbrace 1,2,...,23 \right\rbrace$, $\alpha(3,n+24) = 1$.\\
 Finally, $\alpha(3,n+24) = \alpha(3,n)$, for all $n \in \mathbb{N}^*$, and the homology is of isomorphism type
\begin{align*}
H_1(X_{3,n}) \cong \left\lbrace \begin{matrix}
(\Z/2\Z)^2 \oplus \Z^2 & \text{for} & 24 \nmid n \\
(\Z/2\Z)^2 \oplus (\Z/7\Z)^2 \oplus \Z^2 & \text{for} & 24 \mid n.
\end{matrix} \right.
\end{align*}

With the same procedure as in the proof of the case $m=3$, we can show for arbitrary odd $m$ that
\begin{align*}
p_i \mid \beta(m,n) \ \text{if} \ n = n(p_i)\cdot k \text{ with } k\in \mathbb{N}^*, \ \text{otherwise} \ p_i \nmid \beta(m,n).
\end{align*}
and
\begin{align*}
p_i^{t_i} \mid \beta(m,n) \ \text{if} \ n = p_i^{t_i-1}n(p_i)\cdot k \text{ with } k\in \mathbb{N}^*, \ \text{otherwise} \ p_i^{t_i} \nmid \beta(m,n),
\end{align*}
with $t_i \in \mathbb{N^*}, t_i \ge 2$.\\
From now on, we use the pair $(q,r(q))$ to mean that $q \mid \beta(m,n)$ if $n=r(q)\cdot k$, otherwise $q \nmid \beta(m,n)$.\\
First, we have the set $A$ that contains the pairs
\begin{align*}
(p_i^{t_i}, p_i^{t_i-1}\cdot n(p_i)), \quad i =1,2,...,k; 1 \le t_i \le s_i.
\end{align*}
From these pairs, combining with the properties of prime numbers, we obtain a set of pairs
\begin{align*}
(lcm(B), lcm(C))
\end{align*}
where $B$ is any non empty subset of $A$, and $C = \left\lbrace r(q)| q \in B \right\rbrace$.\\
Sorting these pairs in ascending order of $r(q)$, we may see some pairs
\begin{align*}
(q_1,r(q_1)) \ \text{and} \ (q_2,r(q_2)) \ \text{with} \ r(q_1) = r(q_2).
\end{align*}
Delete these two pairs $(q_1,r(q_1))$ and $(q_2,r(q_2))$ in the list of pairs and replace them by the pair  \\
$(lcm(q_1,q_2),r(q_1))$, because the latter one also satisfies our property $( lcm(q_1,q_2) \mid \beta(m,n)\text{ if }n=r(q_1)\cdot k$, otherwise
$lcm(q_1,q_2) \nmid \beta(m,n)$.
We now sort these pairs in order \underline{strictly} ascending of $r(q)$ and obtain a list of pairs
\begin{align*}
(q_1,r_1), (q_2,r_2), ..., (q_t,r_t)
\end{align*}
for some $t$.
Now,
\begin{itemize}
\item If $n = kr_t$, $\beta(m,n) = q_t$.
\item If $n \neq kr_t$ and $n = kr_{t-1}$, $\beta(m,n) = q_{t-1}$.
\item ...
\item If $n \neq kr_t$,..., $n \neq kr_2$ and $n = kr_1$, then $\beta(m,n) = q_1$.
\item If $n \neq kr_t$,..., $n \neq kr_2$, $n \neq kr_1$, then $\beta(m,n) =1$.
\end{itemize}

Since prime numbers are coprime, $\beta(m,n)$ can be calculated for each pair $(m, n)$ by the formula stated in the theorem, and the period of $\beta(m,n)$ in $n$ is $$r_t = lcm(p_1^{s_1-1}n(p_1),p_2^{s_2-1}n(p_2),...,p_k^{s_k-1}n(p_k)).$$
 \end{proof}

An explicit description of $\beta(m, n)$ at small values of $m$ is given in Table~\ref{explicit-description} above.

\end{document}